\documentclass[12pt]{article}
\usepackage{amssymb,amsmath,indentfirst}
\usepackage[T2A]{fontenc}
\usepackage{srcltx}
\usepackage{graphicx}
\usepackage{xypic}
\newcounter{theorem}
\newcommand{\theorem}{\par\refstepcounter{theorem}%
{\bf Theorem \arabic{theorem}.} }
\input amssym.def
\input amssym.tex

\def\codim{\mathop{\rm codim}}

\title{
{\bf Topological properties of caustics in five-dimensional spaces}
\author{Vyacheslav D. Sedykh}}

\begin{document}

\date{}
\maketitle

\begin{abstract}
We give a list of universal linear relations between the Euler characteristics of manifolds consisting of multisingularities of a generic Lagrangian map into a five-dimensional space. From these relations it follows, in particular, that the numbers $D_5A_2, A_4A_3, A_4A_2^2$ of isolated self-intersection points of the corresponding types on any generic compact four-dimensional caustic are even. The numbers $D_4^+A_3+D_4^-A_3+E_6$, $D_4^+A_2^2+D_4^-A_2^2+\frac12A_4A_3$ are even as well.

\medskip
{\bf Key words:} Lagrangian map, caustic, $ADE$ singularities, multisingularities, adjacency index, Euler characteristic.
\end{abstract}

\section{Introduction}

Caustics are the sets of critical values of so-called Lagrangian maps (see \cite{Arn96}, \cite{Sed21}). Examples of caustics are light caustics, evolutes of plane curves, focal sets of hypersurfaces and other envelopes of systems of rays of different nature.

A generic caustic is a singular hypersurface. The singular points of this hypersurface are described by Arnold's theorem on Lagrangian singularities. Namely germs of a generic Lagrangian map $f:L\rightarrow V$ of a smooth manifold $L$ into a smooth manifold $V$ of the same dimension $n\leq5$ (both without boundary) are Lagrangian equivalent to the germs at the origin of the map
$$
\mathbb{R}^n\to \mathbb{R}^n,\quad
(t,q)\mapsto \left(
-\frac{\partial S(t,q)}{\partial t},q\right),\quad t=(t_1,\dots,t_k),\quad q=(q_{k+1},\dots,q_n)
$$
defined by the function $S=S(t,q)$ of one of the following types corresponding to positive integers $\mu\leq n+1$:
\begin{center}
\begin{tabular}{rll}
$A_{\mu}^\pm:$ & $S=\pm t_1^{\mu+1}+q_{\mu-1}t_1^{\mu-1}+...+q_2t_1^2$, & $\mu\geq1$;\medskip \\
$D_{\mu}^{\pm}:$ & $S=t_1^2t_2\pm t_2^{\mu-1}+q_{\mu-1}t_2^{\mu-2}+...+q_3t_2^2$, & $\mu\geq4$;\medskip \\
$E_6^\pm:$ & $S=t_1^3\pm t_2^4+q_5t_1t_2^2+q_4t_1t_2+q_3t_2^2$, & $\mu=6$.\medskip \\
\end{tabular}
\end{center}

The equivalence class of a Lagrangian map germ at a critical point with respect to Lagrangian equivalence is called a (Lagrangian) singularity. The type of the function $S$ determines the type of a Lagrangian map germ and its singularity. The number $\mu-1$ is called the codimension of a singularity. If $\mu$ is even or $\mu=1$, then germs of types $A_{\mu}^{+}$ and $A_{\mu}^{-}$ are Lagrangian equivalent (their types are denoted by $A_{\mu}$). In other cases, the germs of the types listed above are pairwise Lagrangian non-equivalent.

\medskip
{\bf Remark}. We consider Lagrangian maps in the Whitney $C^{\infty}$-topology. Lagrangian singularities of types $A_{\mu}^\pm,D_{\mu}^\pm,E_6^\pm$ are simple (that is, have zero modality) and stable. For $n>5$ there are singularities that have functional moduli (of $n$ variables) and cannot be removed by a small (Lagrangian) deformation of the mapping.

\medskip
Let $y$ be an arbitrary point of the target space $V$ of a generic proper Lagrangian map $f$. We consider the unordered set of the
symbols from Arnold's theorem that are the types of germs of $f$ at the preimages of $y$. The formal commutative product $\mathcal{A}$ of these symbols is called the type of a multisingularity of $f$ at the point $y$ (or the type of a monosingularity if $y$ has only one preimage). If $f^{-1}(y)=\emptyset$, then $\mathcal{A}=\mathbf{1}$. The types of multisingularities belong to the free Abelian multiplicative semigroup $\mathbb{S}^+$ with the unity $\bf1$ and generators
$$
A_1, A_{2}, A_{4}, A_{6}, A_3^\pm, A_5^\pm, D_4^\pm, D_5^\pm, D_6^\pm, E_6^\pm.
$$

The set $\mathcal{A}_f$ of points $y\in V$ at which $f$ has a multisingularity of type $\mathcal{A}\in\mathbb{S}^+$ is a smooth submanifold of the space $V$. It is called the manifold of multisingularities of type $\mathcal{A}$. The codimension of $\mathcal{A}_f$ in $V$ is equal to the sum of the codimensions of singularities of $f$ at all preimages of an arbitrary point $y\in \mathcal{A}_f$. This sum is called the codimension of a multisingularity of type $\mathcal{A}$ and is denoted by $\mathrm{codim}\,\mathcal{A}$. If $\mathcal{A}_f\neq\emptyset$, then $\mathrm{codim}\,\mathcal{A}\leq n$.

The type of a multisingularity of a Lagrangian map $f$ at a critical value $y$ determines the type of its caustic germ at this point. Germs of a caustic are diffeomorphic if and only if their types either coincide or differ in the number of factors $A_1$ and the signs in the superscript of factors of the forms $A_{2k+1}^\pm,D_{2k+1}^\pm,E_6^\pm$. Therefore the factors $A_1$ and the superscript of the symbols of the mentioned forms are often not written if we are talking about the types of caustic germs and its singular points (but not about the type of Lagrangian map multisingularities). The set of points $y$ at which the caustic of $f$ has a singularity of a given type is a smooth submanifold of $V$. It is the disjoint union of manifolds $\mathcal{A}_f$ such that $\mathcal{A}\in\mathbb{S}^+$ determines the type of the caustic singularity under consideration.

It is well-known (R. Thom -- V. I. Arnold) that the number of singularities of type $A_4$ (swallowtail) and the total number of singularities of types $D_4^+$ (purse) and $D_4^-$ (pyramid) on a generic compact caustic in a three-dimensional space are even. V. A. Vassiliev proved the parity of the following numbers (see \cite{Vas}): the number of type $D_5$ singularities on a generic compact caustic in a four-dimensional space; the total number of types $D_6^\pm$ singularities and the total number of types $A_{6}$ and $E_6$ singularities on a generic compact caustic in a five-dimensional space. 
We proved in \cite[Corollary 16.3]{Sed15} that
the number of type $A_4A_2$ singularities and the total number of types $D_4^\pm A_2$ singularities on a generic compact caustic in a four-dimensional space are even as well. Recently we found new coexistence conditions for isolated singularities of caustics in five-dimensional spaces.

Let the numbers of isolated caustic singularities of types $E_6,D_5A_2,D_4^+A_3,D_4^-A_3,D_4^+A_2^2$, $D_4^-A_2^2,A_4A_3,A_4A_2^2$ be denoted by the same symbols respectively. Then the numbers
\begin{equation}
D_5A_2,\quad A_4A_3,\quad A_4A_2^2,\quad D_4^+A_3+D_4^-A_3+E_6,\quad D_4^+A_2^2+D_4^-A_2^2+\frac12A_4A_3
\label{list-points}
\end{equation}
are even on any generic compact caustic in a five-dimensional space.

This statement is a corollary of Theorem \ref{th-mod}, where two more new congruences modulo $2$ between the Euler characteristics of even-dimensional manifolds of caustic singularities are given. Theorem \ref{th-mod} follows from the universal linear relations between the Euler characteristics of manifolds of caustic singularities listed in Theorem \ref{relat-ca}. The universality of the relations means that their coefficients (rational numbers) do not depend on the mapping.

Finally, Theorem \ref{relat-ca} follows from Theorem \ref{th-relat}, where the list of universal linear relations between the Euler characteristics of manifolds of Lagrangian multisingularities is given. This list is the main result of the paper. We obtained it by computer calculations using the results of the papers \cite{Sed15} and \cite{Sed23} on adjacency indices of Lagrangian multisingularities. The algorithm of the calculations is described in Section \ref{calc}.

\section{The main result}

Suppose that a generic proper Lagrangian map $f:L\rightarrow V$ has a multisingularity of type $X$ at a point $y\in V$ where $\mathrm{codim}\,X=c$. Fix a neighbourhood $U$ of the origin $0$ in $\mathbb{R}^{c}$ and consider a smooth embedding $h:U\to V$ such that $h(0)=y$ and the submanifold $h(U)\subset V$ is transversal to the manifold $X_f$ at $y$. Let $B_{\varepsilon}\subset \mathbb{R}^{c}$ be the open $c$-dimensional ball of radius $\varepsilon>0$ centred at $0$. Then there is a positive number $\varepsilon_0=\varepsilon_0(f,y,h)$ such that for all $\mathcal{A}\in\mathbb{S}^+$ and $\varepsilon<\varepsilon_0$ the set $h(B_{\varepsilon})\cap \mathcal{A}_f$ is a smooth manifold and the equivalence class of this manifold under diffeomorphisms depends only on $\mathcal{A}$ and $X$. We denote this manifold by $\Xi_{\mathcal{A}}(X)$. Its Euler characteristic (the alternating sum of the Betti numbers of the homology groups with compact supports) is denoted by $J_{\mathcal{A}}(X)$.

We say that a multisingularity of type $X$ is adjacent to a multisingularity of type $\mathcal{A}$ if $\mathcal{A}\neq X$ and $\Xi_\mathcal{A}(X)\neq\emptyset$. The number $J_{\mathcal{A}}(X)$ in this case is called the adjacency index of a multisingularity of type $X$ to a multisingularity of type $\mathcal{A}$. The adjacency indices of monosingularities are given in the papers \cite{Sed15} and \cite{Sed23}. The adjacency indices of multisingularities are calculated using \cite[Theorem 5.1]{Sed15}.

Let us assume that the closure of the manifold $\mathcal{A}_f$ is compact. By $\chi_f(\mathcal{A})$ denote the Euler characteristic of $\mathcal{A}_f$. Then
\begin{equation}
\chi_f(\mathcal{A})=\frac12\sum_{X\,\in\,\mathbb{S}^+\setminus\{\mathcal{A}\}}
(-1)^{n-\codim X}J_\mathcal{A}(X)\,\chi_f(X)
\label{relat1}
\end{equation}
if $\codim\mathcal{A}\equiv n-1\,(\mathrm{mod}\,\,2)$. This is the formula (14.1) from the paper \cite{Sed15}. It follows from the fact that the Euler characteristic of an odd-dimensional compact manifold with boundary is half the Euler characteristic of the boundary.

The formula (\ref{relat1}) defines a homogeneous system of linear equations between the Euler characteristics of manifolds of multisingularities of the map $f$. This system can be easily resolved with respect to the Euler characteristics of odd-dimensional manifolds of multisingularities. The corresponding formulas for Lagrangian maps into three- and four-dimensional spaces are given in \cite[Theorems 14.1 and 14.3]{Sed15}. The answer for five-dimensional spaces is given below.

\medskip
{\bf Remark}. We adopt the following conventions: generator $X_\mu^\delta$ of the semigroup $\mathbb{S}^+$ is $X_\mu^+$ if $\delta=+1$ and $X_\mu^-$ if $\delta=-1$; the number $\chi_f(\mathcal{A}A_1^k)$ is equal to zero if $\mathcal{A}$ does not contain generators $A_1$ and $k<0$.

\medskip
Let $L$ and $V$ be smooth five-dimensional manifolds. To simplify formulas below we denote the Euler characteristic of the manifold $\mathcal{A}_f$ by $\mathcal{A}$.

\medskip
\theorem\label{th-relat} {\it The formulas for the Euler characteristics of manifolds $\mathcal{A}_f$ given in Table $\ref{relations4}$ and the first formula in Table $\ref{relations2}$ are valid for $\delta=\pm1$, all non-negative integer $k$ and any generic proper Lagrangian map $f:L\rightarrow V$ such that the set of singular points of its caustic is compact. The second formula in Table $\ref{relations2}$ {\rm(}for the Euler characteristic of the manifold $(A_1^k)_f${\rm)} is valid if $L$ is compact and $k>0$. This formula is valid for $k=0$ if $V$ is compact as well.}

\medskip
The list of formulas in Tables \ref{relations4} and \ref{relations2} was obtained using a computer. In Section \ref{calc} we describe a quick way to obtain the system of equations, which is equivalent to the system of relations (\ref{relat1}).

\medskip
{\bf Remark}. Apparently Theorem \ref{th-relat} describes a complete system (in the sense of \cite{Sed07}) of universal linear relations with real coefficients between the Euler characteristics of manifolds of multisingularities of Lagrangian maps into five-dimensional spaces.

\medskip
Now let $\mathcal{A}\in\mathbb{S}^+$ do not contain the factors $A_1$ and $\mathcal{A}_f^{ca}$ be the disjoint union of manifolds $(\mathcal{A}A_1^k)_f$ over all non-negative integer $k$. To simplify formulas below we also denote the Euler characteristic $\chi_f(\mathcal{A}^{ca})$ of the manifold $\mathcal{A}_f^{ca}$ by $\mathcal{A}$.

\medskip
\theorem\label{relat-ca} {\it The formulas for the Euler characteristics of manifolds $\mathcal{A}_f^{ca}$ given in Table $\ref{relations-ca}$ are valid for $\delta=\pm1$ and any generic proper Lagrangian map $f:L\rightarrow V$ such that the set of singular points of its caustic is compact. The formula in Table $\ref{relations-ca2}$ for the Euler characteristic of the complement ${\bf 1}_f^{ca}$ to the caustic of the map $f$ is valid if $L$ and $V$ are compact.}

\medskip
This statement can be obtained from Theorem \ref{th-relat} using simple calculations. In particular, Theorem \ref{relat-ca} implies:

\medskip
\theorem\label{th-mod} {\it The congruences modulo $2$ between the Euler characteristics of manifolds $\mathcal{A}_f^{ca}$ given in the following list 
\begin{table}[h]
$$
\begin{array}{rl}
&D_6^++D_6^-\equiv0,\quad A_6\equiv E_6^++E_6^-,\quad D_5^+A_2+D_5^-A_2\equiv0,\quad A_4A_2^2\equiv0\\
&\\
&D_4^+A_3^\delta+D_4^-A_3^\delta\equiv E_6^\delta+A_4A_3^-,\quad D_4^+A_2^2+D_4^-A_2^2\equiv A_4A_3^-\\
&\\
&A_4A_3^++A_4A_3^-\equiv0,\quad D_4^++D_4^-\equiv A_4A_3^-, \quad A_4\equiv A_5^+A_2+A_5^-A_2\\
\end{array}
$$
\end{table}

\noindent
are valid for $\delta=\pm1$ and any generic proper Lagrangian map $f:L\rightarrow V$ such that the set of singular points of its caustic is compact.}

\medskip
{\bf Remark}. The first two congruences imply the mentioned above Vassiliev's results on singularities of types $D_6^\pm$ and $A_6,E_6$ on caustics in five-dimensional spaces.

\medskip
Let us transform the formulas in Tables $\ref{relations-ca}$ and $\ref{relations-ca2}$ using the following notations:
$$
A_3^+X+A_3^-X=A_3X,\quad A_5^+X+A_5^-X=A_5X,\quad D_5^+X+D_5^-X=D_5X,
$$
$$
E_6^+X+E_6^-X=E_6X,\quad (A_3^+)^2X+(A_3^-)^2X+A_3^+A_3^-X=A_3^2X
$$
for any $X\in\mathbb{S}^+$ that does not contain factors $A_1$. Then we get:

\medskip
\theorem\label{th-ca-0} {\it All formulas for the Euler characteristics of manifolds of caustic singularities given in Table $\ref{relations-ca-0}$, except the last one, are valid for $\delta=\pm1$ and any generic proper Lagrangian map $f:L\rightarrow V$ such that the set of singular points on its caustic is compact. The last formula is valid if $L$ and $V$ are compact.}

\medskip
Theorem \ref{th-ca-0} is less informative than Theorem \ref{relat-ca}. For example, the parity of the last two numbers from the list (\ref{list-points}) does not follow from the formulas given in Table \ref{relations-ca-0}.

\medskip
{\bf Remark}. The formulas in \cite[Theorem 14.1]{Sed15} for Lagrangian maps into three-dimensional spaces follow from Tables \ref{relations4} and \ref{relations2} if to put $\chi_f(\mathcal{A})=0$ for all $\mathcal{A}\in\mathbb{S}^+$ such that $\codim \mathcal{A}>3$. Similarly, the formulas in \cite[Corollary 16.1]{Sed15} follow from Tables $\ref{relations-ca}$ and $\ref{relations-ca2}$.

\section{Calculations}\label{calc}

Let $X_1,\dots,X_N$ be the generators of the semigroup $\mathbb{S}^+$. By $\mathbb{Z}[\mathbb{S}^+]$ denote the semigroup algebra of $\mathbb{S}^+$. This is the algebra of polynomials of $X_1,\dots,X_N$ with integer coefficients. The types of multisingularities of generic proper Lagrangian map $f:L\rightarrow V$ are monomials $X=X_1^{k_1}\dots,X_N^{k_N}$ in this polynomial algebra.

The substitution $X_i\rightarrow Y_i,i=1,\dots,N$ determines an isomorphism $\mathbb{Z}[\mathbb{S}^+]\rightarrow\mathbb{Z}[\mathbb{S}^+]$. Consider an additive homomorphism
$$
J:\mathbb{Z}[\mathbb{S}^+]\rightarrow\mathbb{Z}[\mathbb{S}^+]
$$
given by the following formula on monomials:
$$
J(X)=\sum\limits_{\mathcal{A}\,\in\,\mathbb{S}^+}(-1)^{\codim \mathcal{A}}\,J_{\mathcal{A}}(X)\mathcal{A}\,,
$$
where $\mathcal{A}=Y_1^{l_1}\dots,Y_N^{l_N}$. From \cite[Theorem 5.1]{Sed15}, it follows that
$$
J(\mathcal{B}\,\mathcal{C})=J(\mathcal{B})J(\mathcal{C})
$$
for all $\mathcal{B},\mathcal{C}\in\mathbb{Z}[\mathbb{S}^+]$. Hence, $J$ is a ring homomorphism. The action of $J$ on the generators of the semigroup $\mathbb{S}^+$ is determined by the formulas in Table \ref{hom-J}, where $\delta=\pm1$. These formulas follow from \cite[Examples 4.2 -- 4.4, Corollaries 7.4, 7.5, 7.11 -- 7.13]{Sed15} and \cite[Corollary]{Sed23}.

Let
$$
\Lambda_n=\sum\limits_{X\,\in\,\mathbb{S}_n^+}(-1)^{\codim X}\,J(X)X,
$$
where $\mathbb{S}_n^+$ is the set of $X\in\mathbb{S}^+$ such that $\codim X\leq n$.
On the one side $\Lambda_n$ is a polynomial of $Y_1,\dots,Y_N$ with coefficients in $\mathbb{Z}[\mathbb{S}^+]$ as polynomials of $X_1,\dots,X_N$:
$$
\Lambda_n=\sum\limits_{\mathcal{A}\,\in\,\mathbb{S}^+}(-1)^{\codim \mathcal{A}}\left(\sum\limits_{X\,\in\,\mathbb{S}_n^+}(-1)^{\codim X}J_{\mathcal{A}}(X)X\right)\mathcal{A}.
$$
The formula (\ref{relat1}) implies
\begin{equation}
\sum\limits_{X\,\in\,\mathbb{S}_n^+}(-1)^{\codim X}J_{\mathcal{A}}(X)\chi_f(X)=(-1)^n\chi_f(\mathcal{A})
\label{relat2}
\end{equation}
for each $\mathcal{A}\in \mathbb{S}_{n-1}^+$ such that $\codim \mathcal{A}\equiv n-1\,(\mathrm{mod}\,\,2)$.

On the other side $\Lambda_n$ is a polynomial of $X_1,\dots,X_N$ with coefficients in $\mathbb{Z}[\mathbb{S}^+]$ as polynomials of $Y_1,\dots,Y_N$. Namely it is equal to the sum
$$
\sum x_1^{k_1}\cdots x_N^{k_N}
$$
over all non-negative integer $k_1,\dots,k_N$ such that $\sum_{i=1}^N k_i\codim X_i\leq n$, where
\begin{equation}
x_i=(-1)^{\codim X_i}\left(\sum\limits_{\mathcal{A}\,\in\,\mathbb{S}^+}(-1)^{\codim \mathcal{A}}\,J_{\mathcal{A}}(X_i)\mathcal{A}\right)X_i.
\label{subs}
\end{equation}
This sum can be obtained from the partial sum of the formal power series
$$
\sum_{k=0}^{+\infty}x_1^k\cdots \sum_{k=0}^{+\infty}x_N^k
$$
of $x_1,\dots,x_N$ by substitution (\ref{subs}).

Now using a computer, we can rewrite the last expression for $\Lambda_n$ as a polynomial of $Y_1,\dots,Y_N$ with coefficients of $X_1,\dots,X_N$. This immediately leads to the system of equations (\ref{relat2}), which is equivalent to the system of equations (\ref{relat1}).

\medskip

\begin{flushleft}

Gubkin University, Moscow, Russia

\medskip

vdsedykh@gmail.com

sedykh@mccme.ru

\end{flushleft}

\begin{table}[p]
$$
\begin{array}{|rl|}
\hline
&\\
D_5^\delta A_1^k=&E_6^\delta A_1^{k-1}+\frac{1}{2}\left(D_6^+A_1^{k-1}+D_6^-A_1^{k-1}+D_5^\delta A_2A_1^{k}+D_5^\delta A_2A_1^{k-2}\right)\\
&\\
A_5^\delta A_1^k=&D_6^-A_1^{k-1}+\frac{1}{2}\left(E_6^+A_1^{k-1}+E_6^-A_1^{k-1}+A_6A_1^{k-1}+A_5^\delta A_2A_1^{k}+A_5^\delta A_2A_1^{k-2}\right)\\
&\\
D_4^\delta A_2A_1^k=&\frac{\delta+1}{2}\left(E_6^+A_1^{k}+E_6^-A_1^{k}\right)+D_6^\delta A_1^{k}+\frac{1}{2}\left(D_5^+A_2A_1^{k-1}+D_5^-A_2A_1^{k-1}\right)\\
&+D_4^\delta A_3^+A_1^{k-1}+D_4^\delta A_3^-A_1^{k-1}+D_4^\delta A_2^2A_1^{k}+D_4^\delta A_2^2A_1^{k-2}\\
&\\
A_4A_2A_1^k=&2\left(E_6^+A_1^{k}+E_6^-A_1^{k}+D_6^+A_1^{k}\right)+A_6A_1^{k}+D_5^+A_2A_1^{k-1}+D_5^-A_2A_1^{k-1}\\
&+A_5^+A_2A_1^{k-1}+A_5^-A_2A_1^{k-1}+A_4A_3^+A_1^{k-1}+A_4A_3^-A_1^{k-1}\\
&+A_4A_2^2A_1^{k}+A_4A_2^2A_1^{k-2}\\
&\\
(A_3^\delta)^2A_1^k=&\frac12\left(E_6^\delta A_1^{k}+D_6^+A_1^{k}+D_6^-A_1^{k}+D_4^+A_3^\delta A_1^{k-1}+3D_4^-A_3^\delta A_1^{k-1}\right.\\
&\left.+A_4A_3^\delta A_1^{k-1}+(A_3^\delta)^2A_2A_1^{k}+(A_3^\delta)^2A_2A_1^{k-2}\right)\\
&\\
A_3^+A_3^-A_1^k=&D_6^-A_1^{k}+\frac12\left(D_4^+A_3^+A_1^{k-1}+D_4^+A_3^-A_1^{k-1}+3D_4^-A_3^+A_1^{k-1}+3D_4^-A_3^-A_1^{k-1}\right.\\
&\left.+A_6A_1^{k}+A_4A_3^+A_1^{k-1}+A_4A_3^-A_1^{k-1}+A_3^+A_3^-A_2A_1^{k}+A_3^+A_3^-A_2A_1^{k-2}\right)\\
&\\
A_3^\delta A_2^2A_1^k=&D_5^+A_2A_1^{k}+D_5^-A_2A_1^{k}+A_5^\delta A_2A_1^{k}+D_4^+A_3^\delta A_1^{k}\\
&+\frac12\left(D_4^+A_2^2A_1^{k-1}+3D_4^-A_2^2A_1^{k-1}+A_4A_3^\delta A_1^{k}+A_4A_2^2A_1^{k-1}\right.\\
&\left.+3A_3^\delta A_2^3A_1^{k}+3A_3^\delta A_2^3A_1^{k-2}\right)+2(A_3^\delta)^2A_2A_1^{k-1}+A_3^+A_3^-A_2A_1^{k-1}\\
&\\
A_2^4A_1^k=&D_4^+A_2^2A_1^{k}+A_3^+A_2^3A_1^{k-1}+A_3^-A_2^3A_1^{k-1}+\frac12\left(A_4A_2^2A_1^{k}+5\sum_{i=0}^1A_2^5A_1^{k-2i}\right)\\
&\\
A_3^{\delta}A_1^k=&\frac{1}{2}\left(D_4^+A_1^{k-1}+3D_4^-A_1^{k-1}+A_4A_1^{k-1}+A_3^{\delta}A_2A_1^k+A_3^{\delta}A_2A_1^{k-2})\right.\\
&-\frac{1}{4}\left(7E_6^{\delta}A_1^{k-1}+3E_6^{-\delta}A_1^{k-1}+7E_6^+A_1^{k-3}+7E_6^-A_1^{k-3}\right.\\
&+4D_6^+A_1^{k-1}+6D_6^+A_1^{k-3}+10D_6^-A_1^{k-1}+8D_6^-A_1^{k-3}\\
&+D_5^+A_2A_1^k+5D_5^+A_2A_1^{k-2}+4D_5^+A_2A_1^{k-4}+2A_6A_1^{k-1}+2A_6A_1^{k-3}\\
&+D_5^{-}A_2A_1^k+5D_5^{-}A_2A_1^{k-2}+4D_5^{-}A_2A_1^{k-4}+\sum_{i=0}^3{3\choose i}A_3^{\delta}A_2^3A_1^{k-2i}\\
&+A_5^{\delta}A_2A_1^k+3A_5^{\delta}A_2A_1^{k-2}+2A_5^{\delta}A_2A_1^{k-4}+A_5^{-\delta }A_2A_1^{k-2}+A_5^{-\delta}A_2A_1^{k-4}\\
&+\sum_{i=1}^2\left(D_4^+A_3^{-\delta}A_1^{k-2i}+3D_4^-A_3^{-\delta}A_1^{k-2i}\right)+9D_4^-A_3^{\delta}A_1^{k-2}+5D_4^-A_3^{\delta}A_1^{k-4}\\
&+3\sum_{i=1}^2D_4^+A_3^{\delta}A_1^{k-2i}+\sum_{i=0}^2{2\choose i}\left(D_4^+A_2^2A_1^{k-2i-1}+3D_4^-A_2^2A_1^{k-2i-1}\right)\\
&-A_4A_3^{\delta}A_1^k+3A_4A_3^{\delta}A_1^{k-2}+2A_4A_3^{\delta}A_1^{k-4}+A_4A_3^{-\delta}A_1^{k-2}+A_4A_3^{-\delta}A_1^{k-4}\\
&\left.+\sum_{i=0}^2{2\choose i}\left(A_4A_2^2A_1^{k-2i-1}+2(A_3^{\delta})^2A_2A_1^{k-2i-1}+A_3^+A_3^-A_2A_1^{k-2i-1}\right)\right)\\
&\\
\hline
\end{array}
$$
\caption{Relations between Euler characteristics $\chi_f(\mathcal{A})$ (part $1$).}
\label{relations4}
\end{table}

\begin{table}[p]
$$
\begin{array}{|rl|}
\hline
&\\
A_2^2A_1^k=&D_4^+A_1^k+A_3^+A_2A_1^{k-1}+A_3^-A_2A_1^{k-1}+\frac{1}{2}\left(A_4A_1^k+3A_2^3A_1^k+3A_2^3A_1^{k-2}\right)\\
&-\frac{1}{4}\left(10E_6^+A_1^k+18E_6^+A_1^{k-2}+10E_6^-A_1^k+18E_6^-A_1^{k-2}+18D_6^-A_1^{k-2}\right.\\
&+10D_6^+A_1^k+16D_6^+A_1^{k-2}+16\sum_{i=0}^1\left(D_5^+A_2A_1^{k-2i-1}+D_5^-A_2A_1^{k-2i-1}\right)\\
&+2A_6A_1^k+6A_6A_1^{k-2}+7\sum_{i=0}^1\left(A_5^+A_2A_1^{k-2i-1}+A_5^-A_2A_1^{k-2i-1}\right)\\
&+12\left(D_4^-A_3^+A_1^{k-3}+D_4^-A_3^-A_1^{k-3}\right)+18D_4^-A_2^2A_1^{k-2}+14D_4^-A_2^2A_1^{k-4}\\
&+2\sum_{i=0}^1(3+2i)\left(D_4^+A_3^+A_1^{k-2i-1}+D_4^+A_3^-A_1^{k-2i-1}\right)+2A_4A_2^2A_1^{k}\\
&+6D_4^+A_2^2A_1^k+18D_4^+A_2^2A_1^{k-2}+12D_4^+A_2^2A_1^{k-4}+12A_4A_2^2A_1^{k-2}\\
&+\sum_{i=0}^1(3+4i)\left(A_4A_3^+A_1^{k-2i-1}+A_4A_3^-A_1^{k-2i-1}\right)+8A_4A_2^2A_1^{k-4}\\
&+8\sum_{i=1}^2\left((A_3^+)^2A_2A_1^{k-2i}+A_3^+A_3^-A_2A_1^{k-2i}+(A_3^-)^2A_2A_1^{k-2i}\right)\\
&\left.+9\sum_{i=0}^2{2\choose i}\left(A_3^+A_2^3A_1^{k-2i-1}+A_3^-A_2^3A_1^{k-2i-1}\right)
+10\sum_{i=0}^3{3\choose i}A_2^5A_1^{k-2i}\right)\\
&\\
A_1^k=&\frac{1}{2}\left(A_2A_1^k+A_2A_1^{k-2}-D_4^{+}A_1^{k-2}-D_4^{+}A_1^{k-4}-3D_4^{-}A_1^{k-2}-D_4^{-}A_1^{k-4}\right)\\
&-\frac{1}{4}\left(\sum_{i=0}^2{2\choose i}\left(A_3^+A_2A_1^{k-2i-1}+A_3^-A_2A_1^{k-2i-1}\right)-A_4A_1^k+2A_4A_1^{k-2}\right.\\
&+A_4A_1^{k-4}+\sum_{i=0}^3{3\choose i}A_2^3A_1^{k-2i}-(E_6^{+}A_1^k+E_6^{-}A_1^k)-9(E_6^{+}A_1^{k-2}\\
&+E_6^{-}A_1^{k-2})-17(E_6^{+}A_1^{k-4}+E_6^{-}A_1^{k-4})-7(E_6^{+}A_1^{k-6}+E_6^{-}A_1^{k-6})\\
&-D_6^{+}A_1^k-8D_6^{+}A_1^{k-2}-15D_6^{+}A_1^{k-4}-6D_6^{+}A_1^{k-6}-A_6A_1^k-2A_6A_1^{k-2}\\
&-5A_6A_1^{k-4}-2A_6A_1^{k-6}-9D_6^{-}A_1^{k-2}-18D_6^{-}A_1^{k-4}-7D_6^{-}A_1^{k-6}\\
&-3(D_5^{+}A_2A_1^{k-1}+D_5^{-}A_2A_1^{k-1})-13(D_5^{+}A_2A_1^{k-3}+D_5^{-}A_2A_1^{k-3})\\
&-15(D_5^{+}A_2A_1^{k-5}+D_5^{-}A_2A_1^{k-5})-5(D_5^{+}A_2A_1^{k-7}+D_5^{-}A_2A_1^{k-7})\\
&-(A_5^{+}A_2A_1^{k-1}+A_5^{-}A_2A_1^{k-1})-5(A_5^{+}A_2A_1^{k-3}+A_5^{-}A_2A_1^{k-3})\\
&-6(A_5^{+}A_2A_1^{k-5}+A_5^{-}A_2A_1^{k-5})-2(A_5^{+}A_2A_1^{k-7}+A_5^{-}A_2A_1^{k-7})\\
&-(D_4^{+}A_3^{+}A_1^{k-1}+D_4^{+}A_3^{-}A_1^{k-1})-5(D_4^{+}A_3^{+}A_1^{k-3}+D_4^{+}A_3^{-}A_1^{k-3})\\
&-7(D_4^{+}A_3^{+}A_1^{k-5}+D_4^{+}A_3^{-}A_1^{k-5})-3(D_4^{+}A_3^{+}A_1^{k-7}+D_4^{+}A_3^{-}A_1^{k-7})\\
&-6(D_4^{-}A_3^{+}A_1^{k-3}+D_4^{-}A_3^{-}A_1^{k-3})-10(D_4^{-}A_3^{+}A_1^{k-5}+D_4^{-}A_3^{-}A_1^{k-5})\\
&-4(D_4^{-}A_3^{+}A_1^{k-7}+D_4^{-}A_3^{-}A_1^{k-7})-6D_4^{+}A_2^2A_1^{k-2}-12D_4^{+}A_2^2A_1^{k-4}\\
&-10D_4^{+}A_2^2A_1^{k-6}-3D_4^{+}A_2^2A_1^{k-8}-6D_4^{-}A_2^2A_1^{k-2}-16D_4^{-}A_2^2A_1^{k-4}\\
&-14D_4^{-}A_2^2A_1^{k-6}-4D_4^{-}A_2^2A_1^{k-8}-3(A_4A_3^{+}A_1^{k-3}+A_4A_3^{-}A_1^{k-3})\\
&-5(A_4A_3^{+}A_1^{k-5}+A_4A_3^{-}A_1^{k-5})-2(A_4A_3^{+}A_1^{k-7}+A_4A_3^{-}A_1^{k-7})\\
&-3A_4A_2^2A_1^{k-2}-8A_4A_2^2A_1^{k-4}-7A_4A_2^2A_1^{k-6}-2A_4A_2^2A_1^{k-8}-D_4^{+}A_2^2A_1^k\\
&-2\sum_{i=1}^4{3\choose i-1}\left((A_3^{+})^2A_2A_1^{k-2i}+A_3^{+}A_3^{-}A_2A_1^{k-2i}+(A_3^{-})^2A_2A_1^{k-2i}\right)\\
&\left.-2\sum_{i=0}^4{4\choose i}\left(A_3^{+}A_2^3A_1^{k-2i-1}+A_3^{-}A_2^3A_1^{k-2i-1}\right)-2\sum_{i=0}^5{5\choose i}A_2^5A_1^{k-2i}\right)\\
&\\
\hline
\end{array}
$$
\caption{Relations between Euler characteristics $\chi_f(\mathcal{A})$ (part $2$).}
\label{relations2}
\end{table}

\begin{table}[p]
$$
\begin{array}{|rl|}
\hline
&\\
D_5^\delta=&E_6^\delta+\frac{1}{2}\left(D_6^++D_6^-\right)+D_5^\delta A_2\\
&\\
A_5^\delta=&D_6^-+\frac{1}{2}\left(E_6^++E_6^-+A_6\right)+A_5^\delta A_2\\
&\\
D_4^\delta A_2=&\frac{\delta+1}{2}\left(E_6^++E_6^-\right)+D_6^\delta+\frac{1}{2}\left(D_5^+A_2+D_5^-A_2\right)\\
&+D_4^\delta A_3^++D_4^\delta A_3^-+2D_4^\delta A_2^2\\
&\\
A_4A_2=&2\left(E_6^++E_6^-+D_6^+\right)+A_6+D_5^+A_2+D_5^-A_2\\
&+A_5^+A_2+A_5^-A_2+A_4A_3^++A_4A_3^-+2A_4A_2^2\\
&\\
(A_3^\delta)^2=&\frac12\left(E_6^\delta+D_6^++D_6^-+D_4^+A_3^\delta+3D_4^-A_3^\delta+A_4A_3^\delta\right)+(A_3^\delta)^2A_2\\
&\\
A_3^+A_3^-=&D_6^-+\frac12\left(D_4^+A_3^++D_4^+A_3^-+3D_4^-A_3^++3D_4^-A_3^-\right.\\
&\left.+A_6+A_4A_3^++A_4A_3^-\right)+A_3^+A_3^-A_2\\
&\\
A_3^\delta A_2^2=&D_5^+A_2+D_5^-A_2+A_5^\delta A_2+D_4^+A_3^\delta+2(A_3^\delta)^2A_2+A_3^+A_3^-A_2\\
&+\frac12\left(D_4^+A_2^2+3D_4^-A_2^2+A_4A_3^\delta+A_4A_2^2\right)+3A_3^\delta A_2^3\\
&\\
A_2^4=&D_4^+A_2^2+A_3^+A_2^3+A_3^-A_2^3+\frac12A_4A_2^2+5A_2^5\\
&\\
A_3^{\delta}=&\frac{1}{2}\left(D_4^++3D_4^-+A_4\right)+A_3^{\delta}A_2\\
&-\frac{1}{2}\left(7E_6^{\delta}+5E_6^{-\delta}+5D_6^++9D_6^-+2A_6\right.\\
&+5D_5^+A_2+5D_5^{-}A_2+3A_5^{\delta}A_2+A_5^{-\delta }A_2\\
&+3D_4^+A_3^{\delta}+7D_4^-A_3^{\delta}+D_4^+A_3^{-\delta}+3D_4^-A_3^{-\delta}\\
&+2D_4^+A_2^2+6D_4^-A_2^2+2A_4A_3^{\delta}+A_4A_3^{-\delta}\\
&\left.+2A_4A_2^2+4(A_3^{\delta})^2A_2+2A_3^+A_3^-A_2+4A_3^{\delta}A_2^3\right)\\
&\\
A_2^2=&D_4^++A_3^+A_2+A_3^-A_2+\frac{1}{2}A_4+3A_2^3\\
&-\frac{1}{2}\left(14E_6^++14E_6^-+13D_6^++9D_6^-+4A_6\right.\\
&+16\left(D_5^+A_2+D_5^-A_2\right)+7\left(A_5^+A_2+A_5^-A_2\right)\\
&+8\left(D_4^+A_3^++D_4^+A_3^-\right)+6\left(D_4^-A_3^++D_4^-A_3^-\right)\\
&+18D_4^+A_2^2+16D_4^-A_2^2+5A_4A_3^++5A_4A_3^-\\
&+11A_4A_2^2+8\left((A_3^+)^2A_2+A_3^+A_3^-A_2+(A_3^-)^2A_2\right)\\
&\left.+18A_3^+A_2^3+18A_3^-A_2^3+40A_2^5\right)\\
&\\
\hline
\end{array}
$$
\caption{Relations between Euler characteristics $\chi_f(\mathcal{A}^{ca})$ (part $1$).}
\label{relations-ca}
\end{table}

\begin{table}[p]
$$
\begin{array}{|rl|}
\hline
&\\
\mathbf{1}=&A_2-D_4^{+}-2D_4^{-}-\frac{1}{2}A_4-A_3^+A_2-A_3^-A_2-2A_2^3\\
&+\frac{1}{2}\left(17E_6^{+}+17E_6^{-}+15D_6^{+}+17D_6^{-}+5A_6\right.\\
&+18D_5^{+}A_2+18D_5^{-}A_2+7A_5^{+}A_2+7A_5^{-}A_2\\
&+8D_4^{+}A_3^{+}+8D_4^{+}A_3^{-}+10D_4^{-}A_3^{+}+10D_4^{-}A_3^{-}\\
&+16D_4^{+}A_2^2+20D_4^{-}A_2^2+5A_4A_3^{+}+5A_4A_3^{-}\\
&+10A_4A_2^2+8\left((A_3^{+})^2A_2+A_3^{+}A_3^{-}A_2+(A_3^{-})^2A_2\right)\\
&\left.+16A_3^{+}A_2^3+16A_3^{-}A_2^3+32A_2^5\right)\\
&\\
\hline
\end{array}
$$
\caption{Relations between Euler characteristics $\chi_f(\mathcal{A}^{ca})$ (part $2$).}
\label{relations-ca2}
\end{table}

\begin{table}[p]
$$
\begin{array}{|rl|}
\hline
&\\
D_5=&E_6+D_6^++D_6^-+D_5A_2\\
&\\
A_5=&2D_6^-+E_6+A_6+A_5A_2\\
&\\
D_4^\delta A_2=&\frac{\delta+1}{2}E_6+D_6^\delta+\frac{1}{2}D_5A_2+D_4^\delta A_3+2D_4^\delta A_2^2\\
&\\
A_4A_2=&2\left(E_6+D_6^+\right)+A_6+D_5A_2+A_5A_2+A_4A_3+2A_4A_2^2\\
&\\
A_3^2=&\frac12\left(E_6+A_6\right)+D_6^++2D_6^-+D_4^+A_3+3D_4^-A_3+A_4A_3+A_3^2A_2\\
&\\
A_3A_2^2=&2D_5A_2+A_5A_2+D_4^+A_3+D_4^+A_2^2+3D_4^-A_2^2\\
&+\frac12A_4A_3+A_4A_2^2+2A_3^2A_2+3A_3A_2^3\\
&\\
A_2^4=&D_4^+A_2^2+A_3A_2^3+\frac12A_4A_2^2+5A_2^5\\
&\\
A_3=&D_4^++3D_4^-+A_4+A_3A_2\\
&-\left(6E_6+5D_6^++9D_6^-+2A_6+5D_5A_2+2A_5A_2\right.\\
&+2D_4^+A_3+5D_4^-A_3+2D_4^+A_2^2+6D_4^-A_2^2\\
&\left.+2A_4A_2^2+2A_3^2A_2+2A_3A_2^3\right)-\frac{3}{2}A_4A_3\\
&\\
A_2^2=&D_4^++A_3A_2+\frac{1}{2}A_4+3A_2^3\\
&-\frac{1}{2}\left(14E_6+13D_6^++9D_6^-+4A_6+16D_5A_2+7A_5A_2\right.\\
&+8D_4^+A_3+6D_4^-A_3+18D_4^+A_2^2+16D_4^-A_2^2+5A_4A_3\\
&\left.+11A_4A_2^2+8A_3^2A_2+18A_3A_2^3+40A_2^5\right)\\
&\\
\mathbf{1}=&A_2-D_4^{+}-2D_4^{-}-\frac{1}{2}A_4-A_3A_2-2A_2^3\\
&+\frac{1}{2}\left(17E_6+15D_6^{+}+17D_6^{-}+5A_6+18D_5A_2+7A_5A_2\right.\\
&+8D_4^{+}A_3+10D_4^{-}A_3+16D_4^{+}A_2^2+20D_4^{-}A_2^2+5A_4A_3\\
&\left.+10A_4A_2^2+8A_3^2A_2+16A_3A_2^3+32A_2^5\right)\\
&\\
\hline
\end{array}
$$
\caption{Relations between Euler characteristics of manifolds of caustic singularities.}
\label{relations-ca-0}
\end{table}

\begin{table}[p]
$$
\begin{array}{|rl|}
\hline
&\\
J(1)&=1\\
&\\
J(A_1)&=A_1\\
&\\
J(A_2)&=1+A_1^2-A_2\\
&\\
J(A_3^\delta)&=A_1+A_1^3-2A_2A_1+A_3^\delta\\
&\\
J(A_4)&=1+A_1^2+A_1^4-A_2-3A_2A_1^2+A_2^2+A_3^+A_1+A_3^-A_1-A_4\\
&\\
J(D_4^\delta)&=\frac{1+\delta}{2}\left(1+2A_1^2-2A_2+2A_2^2\right)+\frac{3-\delta}{2}A_1^4\\
&\quad-(5-\delta)A_2A_1^2+(2-\delta)\left(A_3^+A_1+A_3^-A_1\right)-D_4^{\delta}\\
&\\
J(A_5^\delta)&=A_1+A_1^3+A_1^5-2A_2A_1-4A_2A_1^3+3A_2^2A_1\\
&\quad+A_3^\delta+2A_3^\delta A_1^2+A_3^{-\delta}A_1^2-2A_3^\delta A_2-2A_4A_1+A_5^\delta\\
&\\
J(D_5^\delta)&=A_1+A_1^3+2A_1^5-3A_2A_1-9A_2A_1^3+6A_2^2A_1+A_3^++A_3^-\\
&\quad+4A_3^+A_1^2+4A_3^-A_1^2-2A_3^+A_2-2A_3^-A_2-2A_4A_1-D_4^+A_1-D_4^-A_1+D_5^\delta\\
&\\
J(A_6)&=1+A_1^2+A_1^4+A_1^6-A_2-3A_2A_1^2-5A_2A_1^4+A_2^2+6A_2^2A_1^2\\
&\quad+A_3^+A_1+A_3^-A_1+2A_3^+A_1^3+2A_3^-A_1^3-A_2^3-3A_3^+A_2A_1-3A_3^-A_2A_1\\
&\quad-A_4-3A_4A_1^2+2A_4A_2+A_3^+A_3^-+A_5^+A_1+A_5^-A_1-A_6\\
&\\
J(D_6^\delta)&=\frac{1+\delta}{2}(1+2A_1^2+2A_1^4-2A_2+3A_2^2-4A_2^3-2A_4+4A_4A_2)+\frac{5-\delta}{2}A_1^6\\
&\quad+\frac{33-\delta}{2}A_2^2A_1^2+\frac{5+\delta}{2}(A_3^+A_1+A_3^-A_1)+\frac{13-3\delta}{2}(A_3^+A_1^3+A_3^-A_1^3)\\
&\quad-(5+3\delta)A_2A_1^2-2(7-\delta)A_2A_1^4-(9-2\delta)(A_3^+A_2A_1+A_3^-A_2A_1)\\
&\quad+(A_3^+)^2+(A_3^-)^2+(1-\delta)(A_3^+A_3^-+A_5^+A_1+A_5^-A_1)-(5-\delta)A_4A_1^2\\
&\quad-D_4^\delta-2D_4^{\delta}A_1^2-D_4^{-\delta}A_1^2+2D_4^\delta A_2+D_5^+A_1+D_5^-A_1-D_6^\delta\\
&\\
J(E_6^\delta)&=1+2A_1^2+A_1^4+2A_1^6-2A_2-6A_2A_1^2-12A_2A_1^4+3A_2^2+16A_2^2A_1^2-4A_2^3\\
&\quad+2A_3^\delta A_1+3A_3^{-\delta}A_1+6A_3^\delta A_1^3+5A_3^{-\delta}A_1^3-8A_3^\delta A_2A_1-6A_3^{-\delta}A_2A_1+(A_3^\delta)^2\\
&\quad-2A_4-6A_4A_1^2+4A_4A_2+A_5^+A_1+A_5^-A_1\\
&\quad-D_4^+-2D_4^+A_1^2-D_4^-A_1^2+2D_4^+A_2+2D_5^{\delta}A_1-E_6^{\delta}\\
&\\
\hline
\end{array}
$$
\caption{The action of the homomorphism $J$.}
\label{hom-J}
\end{table}

\end{document}